\documentclass[12pt,leqno]{article} 
\usepackage{amsmath,amsfonts,amscd,ProtocoloTope,url}

\title{\Large\bf  Top dimensional group of the basic intersection 
cohomology for singular riemannian foliations\footnote{To appear in the Bulletin of the Polish Academy of Sciences.}.}

\author{
JosŽ Ignacio Royo Prieto\thanks{Departamento de  Matem‡tica Aplicada. 
Universidad del Pa's Vasco. Alameda de Urquijo s/n.  Bilbao - Espa–a.
Partially supported by the UPV-EHU 
grant 
127.310-E-14790/2002, by a PostGrant from 
the Gobierno Vasco-Eusko Jaurlaritza  and by the MCyT of the Spanish 
Government. {\sl joseignacio.royo@ehu.es}}
\\ {\small Universidad del Pa's Vasco}
\and
Martintxo Saralegi-Aranguren\thanks{
FŽdŽration CNRS Nord-Pas-de-Calais FR 2956.
UPRES-EA 2462 LML. 
FacultŽ Jean Perrin.
UniversitŽ d'Artois.   Rue Jean Souvraz SP 18.   62 307 Lens Cedex - 
France. Partially supported by the UPV-EHU 
grant 
127.310-E-14790/2002. 
{\sl saralegi@euler.univ-artois.fr}}
\\ {\small UniversitŽ d'Artois }
\and   
Robert Wolak\thanks{Instytut Matematyki. Uniwersytet Jagiellonski. Wl. 
Reymonta4, 30 059
Krakow - Poland.  
Partially supported by the KBN grant 5 PO3A 
02125. {\sl robert.wolak@im.uj.edu.pl} }
\\ {\small Uniwersytet Jagiellonski}
}

\begin{document}  
\date{}
\maketitle

\begin{abstract}
    It is known that, for a regular riemannian foliation on a compact 
    manifold, the properties of its basic cohomology (non-vanishing 
    of the top-dimensional group and PoincarŽ duality) and the 
    tautness of the foliation are closely related. If we consider 
    singular riemannian foliations, there is little or no relation between 
    these properties.

    We present an example of a singular isometric flow for which the 
    top dimensional basic cohomology group is non-trivial, but its 
    basic cohomology does not satisfy the PoincarŽ Duality property.
    We recover this property in the basic intersection cohomology. 
    
    It 
    is not by chance that the top dimensional basic intersection 
    cohomology groups of the example are isomorphic to either $0$ or $\R$.
    We prove in this Note that this holds for any singular riemannian 
    foliation of a compact connected manifold.
    As a Corollary, we get that the tautness of the regular 
    stratum of the singular riemannian foliation can be detected by  the basic 
    intersection cohomology.
    \end{abstract}
For regular riemannian foliations on compact connected manifolds there 
is a very close relation between tautness (the existence of a 
bundle-like metric for which all leaves are minimal submanifolds) and 
the properties of the basic cohomology. In fact, it was shown that for 
a regular riemannian foliation $\F$ of codimension $n$ on a compact 
connected manifold $M$ the following conditions are equivalent
(cf. \cite{D,EH,KT,MA,To})
\begin{itemize}
    \item[$(i)$] $\F$ is taut.
    \item[$(ii)$] $\hiru{H}{n}{\mf}\neq 0$.
    \item[$(iii)$] $\hiru{H}{n}{\mf} = \R$.
    \item[$(iv)$] $\hiru{H}{*}{\mf}$ verifies the PoincarŽ Duality (PD)
    property,
    \end{itemize}
    under some assumptions of orientability.
Although the basic cohomology of a singular 
riemannian foliation (SRF for short) on a compact manifold is finite 
dimensional (cf. \cite{W}), the relation between the above conditions 
is not straightforwardly exportable to the singular framework.
As pointed out by the authors in \cite {MW}, the singular nature of a 
SRF on a compact manifold  prevents any global metric on it from 
making all the leaves minimal (see also \cite{RSW1}). 

The example 
presented in this paper shows that if, for flows, we replace the 
condition ``taut'' by ``isometric'' - these two conditions are 
equivalent for regular riemannian flows (see \cite{MS}) - we cannot recover the PD 
property. To recover this property  we have to adapt the basic 
cohomology to the stratification defined by the SRF. With that purpose, 
 we have 
defined  the ``basic intersection cohomology'' (BIC for short) for SRFs 
(see \cite{SW1,SW2}). The calculations for the example show that its
BIC satisfies the PoincarŽ Duality. It has been proven 
in \cite{RP} that the BIC of any singular riemannian flow (isometric 
or not) satisfies this property.

The main part of this Note is dedicated to the proof of the 
equivalence of conditions $(i)$, $(ii)$ and $(iii)$ for SRFs on compact 
connected manifolds for the BIC.
We also prove that the top dimensional BIC groups are isomorphic to 
$0$ or $\R$.

\smallskip
%
%
\medskip

The authors would like to thank the referee of the paper for many 
useful comments which helped improve the paper.

		\bigskip

		In the sequel $M$  is a connected, second countable, Haussdorff,   
		without boundary and smooth 
		(of class $C^\infty$) manifold  of dimension
		$m$. 
		All the maps are considered smooth unless something else is indicated.
If $\F$ is a foliation on $M$ and $V$ is a saturated submanifold of $M$ we shall write 
    $(V,\F)$ the induced foliated manifold and $\F_{V}$ the induced 
 foliation.

     \section{Presentation of the singular riemannian 
     foliations\protect\footnote{For the notions related 
    with riemannian foliations we refer the reader to \cite{Mo,To}
    and for the notions related with singular riemannian foliations we 
    refer the reader to \cite{Mo,Mo0,Mo1,BM}.}}
    We are going to work in the framework of the singular riemannian 
    foliations introduced by Molino.
    
  \prg {\bf The SRF}.   A {\em singular riemannian  foliation} (SRF for short) on a manifold $M$ is 
    a partition $\mathcal{F}$ by connected immersed 
    submanifolds, called {\em leaves}, verifying the following properties:

    \begin{itemize}
    \item[I-] The module of smooth vector fields tangent to the leaves is 
    transitive on each leaf.

    \item[II-] There exists a riemannian metric $\mu$ on $M$, called {\em 
    adapted metric},  such that each geodesic 
    that is perpendicular at one point to a leaf remains perpendicular to every 
    leaf it meets.
    \end{itemize}
    The first condition implies that $(M,\mathcal{F})$ is a singular foliation 
    in the sense of \cite{Su} and \cite{St}. Notice that the 
    restriction of $\mathcal{F}$ to a saturated open subset 
    produces a SRF. 
    
    \bigskip
    
    In the next two subsections we recall some basic properties of 
    SRFs which can be found in      \cite{Mo,BM} or easy corollaries of the 
    properties proved there.

  \prg {\bf Stratifications}.   Classifying the points of
  $M$ following the dimension of the leaves one gets a stratification
  $\SF$ of 
  $M$ whose elements are called {\em strata}. The strata are smooth 
  embedded submanifolds. The restriction of 
  $\mathcal{F}$ to a stratum $S$ is a regular foliation 
  $\mathcal{F}_{S}$. The strata are ordered by: $S_{1} \preceq S_{2} 
  \Leftrightarrow S_{1}\subset \overline{S_{2}}$. 
  
  There are several types of strata. The minimal (resp. 
  maximal) strata are the closed strata (resp. open strata). The open 
  strata are called  {\em regular strata} and the others
  are called {\em singular strata}. We 
  denote by   $\SiF$ the family of singular strata. 
  In the case of SRFs, the singular strata 
  are of codimension greater than 1, so there is just one regular 
  stratum, if the manifold is connected (cf. \cite{Mo}).
  The {\em dimension of the foliation} $\mathcal{F}$ is the dimension of the 
  biggest leaves of $\mathcal{F}$, that is, $\dim \mathcal{F} = \dim 
  \mathcal{F}_{R}$.
  The union of singular strata is the {\em singular part} 
  $\Sigma = M \backslash  R.$
  
  The 
  stratum $S$ is a {\em boundary stratum} 
  if there exists  a stratum $S'$ with $S\preceq S'$ and $\codim_{M} \F = 
  \codim_{S'}\F_{S'}- 1$. The reason for this term  can be
explained by the following example.
  Take $M = \S^{4} $ and $\F$ given by the orbits of 
	 the $ \T^{2}$-action:   $(u,v)\cdot (z_{1},z_{2},t) =   (u \cdot z_{1}, 
	 v \cdot z_{2},t),
	 $ where $ \S^4 = \{ (z_{1},z_{2},t) \in \C ^{2} \times \R
	 \ \big| \ |z_{1}|^{2}  +|z_{2}|^{2} + t^{2}= 1 \}. $ Here, the north pole 
	 $(0,0,1)$,
	the
	 south pole $(0,0,-1)$ and the cylinders $\{ z_{1}\ne 0\}$, $\{z_{2} 
	 \ne 0\}$
	   are the boundary strata 
	 and we have $
	 M/\F = \D = \{ (x,y) \in \R^{2}\ \big| \ x^{2}+y^{2} \leq 1\}$ . 
	 The boundary $\partial ( M/\F )$ is given by $\S^{1}$, the union of 
	 the quotient of the boundary strata.
	 In fact, the link of the maximal boundary strata is a sphere  with the 
one leaf foliation (see for example \cite{SW1} for the notion 
	 of  a link).
  
  The  {\em depth}  of $\SF$, written $ \depth  \SF$, is defined to be the largest 
$i$ for which there 
exists a chain of strata $S_0  \prec S_1 \prec \cdots \prec S_i$. So, 
$ \depth  \SF= 0$ if and only if the foliation 
${\mathcal  F}$ is regular.

    \prg {\bf Tubular neighborhood.} 
    Any stratum $S \in \SF$ is a proper submanifold of the 
    riemannian manifold $(M,\F,\mu)$.
    So, it possesses a tubular neighborhood $(T_S,\tau_S,S)$.
    Recall that associated with this neighborhood there are the following smooth maps:
    \begin{itemize}
	\item[+] The {\em radius map} $\rho_S \colon T_S \to [0,1[$ defined 
	fiberwise:
    $z\mapsto |z|$. Each $t\not= 0$ is a regular value of the $\rho_S$. 
    The pre-image $\rho_S^{-1}(0)$ is $S$.
	 \item[+] The {\em contraction} $H_S \colon T_S \times [0,1] \to T_S$ 
	 defined fiberwise: $(z,r) \mapsto r \cdot z$. The restriction 
	 $(H_S)_t \colon T_S \to T_S$ is an 
	 embedding for each 
	 $t\not= 0$   and $(H_{S})_0 \equiv \tau_S$. 
    \end{itemize}

    \nt These two maps verify $\rho_S(r \cdot u) = r  \rho_S(u)$. 
This tubular neighborhood can be chosen verifying the two following important 
properties:

    \Zati Each  $(\rho_S^{-1}(t),\mathcal{F})$ is a 
	SRF, and

	 \zati Each  $(H_S)_{t} \colon  (T_S,\mathcal{F}) \to (T_S,\mathcal{F}) $ is a
	 foliated map.

	 \medskip

	 \nt 
We shall say that $(T_S,\tau_S,S)$   is a {\em foliated tubular 
neighborhood} of $S$. The existence of foliated tubular neighborhoods 
follows from the 
homothetic transformation Lemma of \cite{Mo}.

The hypersurface $D_S = \rho_S^{-1}(1/2)$ is the {\em 
	 core} of the tubular neighborhood. We have the inequality $\depth \SFDS < \depth 
	 \SFTS$. 

\section{Presentation of the basic intersection cohomology\protect\footnote{For 
the notions related with the basic cohomology we refer the reader to 
\cite{ESH,W}, for the notions related 
    with the basic intersection cohomology we refer the reader to 
    \cite{SW1,SW2} and for the notions related with the 
    intersection cohomology we refer the reader to  \cite{GM,Bry}.}}

    Goresky and MacPherson introduced the intersection cohomology for 
    the study of the singular manifolds.
    This cohomology generalizes 
    the usual deRham cohomology for manifolds and possesses similar 
    properties. Following the same principle, 
    the basic intersection cohomology has been 
    introduced for the study of SRFs generalizing the basic cohomology.
    
    \medskip
    
    We fix for the sequel a manifold 
    $M$ endowed with a SRF $\mathcal{F}$. We write $m =\dim M$ and $n = \codim_M \mathcal{F}$.
    
    \prg {\bf The BIC}.  A {\em perversity} is a map 
    $\per{p} \colon \SiF \to \overline{\Z} = \Z \cup \{ -\infty, \infty\}$. 
     There are several particular 
     perversities,
     \begin{itemize}
     \item[-] the {\em constant perversity} $\per{c}$,  defined by 
     $\per{c}(S) =c$, where  
     $c\in \overline{\Z}$, 
     \item[-] the {\em 
     (basic) top 
     perversity} $\per{t}$, defined by 
     $\per{t}(S) =  n - 
     \codim_{S}    \mathcal{F}_S - 2,
     $
     (cf. 1.2) and
     \item[-] the {\em boundary perversity} $\per{\partial}$, defined by
     $
     \per{\partial} = \min (\per{0},\per{t}).
 	$
     \end{itemize}

    The {\em basic intersection cohomology} (BIC for short) 
    $\lau{\IH}{*}{\per{p}}{\mf}$ is the cohomology of the complex 
    $\lau{\Om}{*}{\per{p}}{\mf}$ 
    of 
    {\em $\per{p}$-intersection basic forms}.
    A $\per{p}$-intersection  basic form 
    is a basic form defined on $R$ possessing a 
    vertical degree $\| \omega \|_{S}$, relatively to a foliated 
    tubular neighborhood $(T_{S},\tau_{S},S)$, lower 
    than $\per{p}(S)$ and this for each singular stratum $S$ (see \cite {SW1,SW2} for the exact definition).
    Recall that $\| \omega \|_{S} \leq i$ when $\omega(v_{0}, \ldots 
    , v_{i}, -) = 0$ for each family $\{v_{0}, \ldots , v_{i}\}$ of 
    vectors tangent to the fibers of $\tau_{S}$.
    
 When $\depth \SF =0$ then we have $\lau{H}{*}{\per{p}}{\mf} = 
 \hiru{H}{*}{\mf}$ for any perversity.
    
%
%
%
%
%
%
    \prg {\bf Mayer-Vietoris}. 
    A  covering $\{ U , V \}$  of $M$ by saturated open 
     subsets possesses a 
     subordinated partition of the unity made up of  basic  
     functions (see Lemma below). 
     	
	For a such covering we have the Mayer-Vietoris short sequence
\begin{equation}
\label{mv}
	0 \to \lau{\Om}{*}{\per{p}}{\mf} \to 
	\lau{\Om}{*}{\per{p}}{U/\mathcal{F}} \oplus 
	\lau{\Om}{*}{\per{p}}{V/\mathcal{F}}  \to
	\lau{\Om}{*}{\per{p}}{(	U \cap V)/\mathcal{F}}  \to 0,
\end{equation}
	where the maps are defined by restriction. The third map is onto since the 
	elements of the partition of the unity are $\per{0}$-basic  
	functions and $\lau{\Om}{*}{\per{0}}{\bullet/\F} \cdot 
	\lau{\Om}{*}{\per{p}}{\bullet/\F} \subset \lau{\Om}{*}{\per{p}}{\bullet/\F} $. Thus,  
	the sequence
	 is exact. This result is not 
	longer true for 
	more general coverings.   
	
	\medskip

	For the existence of the above Mayer-Vietoris sequence   we need the 
	following folk result, well-known for compact Lie group actions and 
	regular riemannian foliations.
	\bL
	\label{sat}
	A covering $\{ U , V \}$  of $M$ by saturated open 
	subsets possesses a 
	subordinated partition of the unity made up of  basic  
	functions. 
	\eL

    \prg {\bf Compact supports.} In this Note we need to work with the BIC 
   with compact supports. The {\em support} of a differential form $\om \in 
    \lau{\Om}{*}{\per{p}}{\mf}$, written $\supp \om$, is the closure in 
$M$ of $\{ x \in M \ | \ \om(x) \ne 0 \}$. 
We denote by  $\lau{\Om}{*}{\per{p},c}{\mf}$ the complex 
of $\per{p}$-intersection basic  forms with compact support. Its cohomology is 
$\lau{\IH}{*}{\per{p},c}{\mf}$. When $M$ is compact, we have 
$\lau{\IH}{*}{\per{p},c}{\mf}= \lau{\IH}{*}{\per{p}}{\mf}$ and 
when $\depth \SF =0$ then we have $\lau{H}{*}{\per{p},c}{\mf} = 
\lau{H}{*}{c}{\mf}$, for any perversity.

  Associated to a saturated open  covering $\{U,V\}$ of $M$ we have
the Mayer-Vietoris short sequence (see Lemma above)
    \begin{equation}
\label{MV}
    0 \to 
  \lau{\Om}{*}{\per{p},c}{(U \cap V)/\mathcal{F}}  
\to 
    \lau{\Om}{*}{\per{p},c}{U/\mathcal{F}} 
\oplus 
    \lau{\Om}{*}{\per{p},c}{V/\mathcal{F}}  
\to
\lau{\Om}{*}{\per{p},c}{\mf}\to 0,
    \end{equation}
    where the map are defined by inclusion. The third map is onto since the 
    elements of the partition of the unity are $\per{0}$-basic  
    functions. Thus,  
    the sequence
     is exact. 

\prg {\bf Example}. 
	       us consider the isometric action
	  $
	  \Phi \colon \R \times \S^{2d+2} \to \S^{2d+2}
	  $
	  given by the formula 
	  $$
	  \Phi (t, (z_{0}, \ldots , z_{d},x)) = (e^{a_{0} \pi it} \cdot z_{0}, \ldots , 
	  e^{a_{d}\pi it}\cdot z_{d} ,x),
	  $$
	  with $(a_{0}, \ldots , a_{d}) \ne (0,\ldots , 0)$. Here, 
	  $
	  \S^{2d+2} = \{ (z_{0}, \ldots , z_{d},x) \in \C^d \times \R \  |  \ |z_{0}|^{2} + \cdots +
	  |z_{d}|^{2}+ x^{2}= 1 \}$.
	  There are two singular strata: the north pole $S_{1}= (0, \ldots , 
	  0,1)$ and the south pole $S_{2} = (0,\ldots,-1)$.
	  the regular stratum is $R = \S^{2d+1} \times ]-1,1[$. Let  $r \colon 
	  R\to \R$ a smooth map, depending  just on 
	  $]1,1[$, 
	  with $r \equiv 0$ on $]0,1/4[ \cup ]3/4,1[$ and $\Int{0}{1} r= 1$.
	  The basic cohomology $\hiru{H}{*}{\S^{2d+2}/\F}$ of
	  the foliation defined by this 
	  flow is the following:

	  \medskip

	  \begin{center}
	  \begin{tabular}{|c|c|c|c|c|c|c|c|c|c|}\hline
	      &&&&&&&&\\
	   $i=0$ & $i=1$ & $i=2$  & $i=3$ & $i=4$ & $i=5$  & $i= \cdots $ & 
	   $i=2d$ & $i=2d+1$    \\[,1cm] \hline
	   &&&&&&&&\\
	    $1$&  0 & 0  & $ [dr
	  \wedge e] $ 
	  &0 & $[dr \wedge e^2] $ & $\cdots$ & 0 &  $ [dr \wedge e^d] $    \\[,1cm]  \hline
	  \end{tabular}
	  \end{center}

	  \medskip

	  \nt where $e \in \lau{\Om}{2}{\per{2}}{\S^{2d+2}/\mathcal F}$
	  is an {\em Euler 
	  form} (cf. \cite{HS}). 
	  These calculations come directly from the equalities:
	  \begin{equation}
	      \label{calculos}
	      \| r \|_{S_{k}} = 0,  
	      \| dr \|_{S_{k}} = -\infty,
	      \| e^j \|_{S_{k}} = 
	      \| re^j\|_{S_{k}} = 
	  2^j, \ 
	  \|dr \wedge e^j\|_{S_{k}} = -\infty \ 
	  \hbox{ and } \ 
	  dr \wedge e^j= d(re^j),
	  \end{equation}
	  for $k=1,2$ and $j \in \{ 1 , \ldots, d\}$. 

	  We notice that the top dimensional basic cohomology group is 
	  isomorphic to $\R$, but this cohomology does not have the PoincarŽ 
	  Duality property  in spite of the fact that the flow is isometric!
	  The classical basic cohomology does not take 
	  into account the stratification $\SF$. 
	  However,  even for the SRF, that basic cohomology is finite 
	  dimensional (cf. \cite{W}).

	  If we consider the BIC of our  example the picture changes. 
	  The following table presents the BIC
	  $\lau{\IH}{*}{\per{p}}{\S^{k =2d+2}/\mathcal{F}}$ for the constant 
	  perversities:

	\medskip
	
	  \begin{center}
	  \begin{tabular}{|c|c|c|c|c|c|c|c|c|c|c|c|c|}\hline
	      &&&&&&&&&&&\\
      i \! = \!& $ 0$ & $ 1$ & $ 2$  & $ 3$ & $ 4$ & $ 5$  &$6$& $7$& $  \cdots $ & $ 
      k \! - \! 2$ & $    k \! - \!  1$
      \\[,1cm] \hline
      &&&&&&&&&&&\\
	   $\per{p} <\per{0}$&0&$[dr]$ & 0&$[e \wedge dr]$ &  0 & $[e^2 \wedge dr]$ 
	   &  0 & $[e^3 \wedge dr]$ &$\cdots$ &0& $[e^d \wedge dr]$
	   \\[,1cm]  \hline
	   &&&&&&&&&&&\\
	   $\per{p}  \! = \!\per{0},\per{1}$ &  $1$ & 0 &  0 &  $[e \wedge dr]$ &  0 &$
		   [e^2 \wedge dr]$&  0 & $[e^2 \wedge dr]$ & $\cdots$  & 0&$[e^d \wedge dr]$ 
		   \\[,1cm]  \hline
		   &&&&&&&&&&&\\
	  $\per{p}  \! = \!\per{2},\per{3}$&0&0 & [e]& 0 &  0 & $[e^2 \wedge dr]$&  0 & $[e^3 \wedge dr]$ &$\cdots$ 
	  &0&$[e^d \wedge dr]$ 
	  \\[,1cm]  \hline
	  &&&&&&&&&&&\\
	  $\per{p}  \! = \!\per{4},\per{5}$ &  $1$ & 0 &  $[e]$ & $0$ &  $[e^{2}]$ & $0$&  0 & $[e^3 \wedge dr]$ 
	  &$\cdots$ &0&$[e^d \wedge dr]$
	  \\[,1cm]  \hline      &&&&&&&&&&&\\
	  $\cdots$ & $\cdots$ & $\cdots$ & $\cdots $& $\cdots$ & $\cdots$& $\cdots$ & $\cdots$ & 
	  $\cdots$  & $\cdots$& $\cdots$& $\cdots$\\[,1cm] 
	  \hline       &&&&&&&&&&&\\
	  $\per{p}  \! = \!\per{k  \! \! - \!   \!  4}, \per{k \! \! - \! \! 
	  3}$ &  $1$ & 0 &  $[e]$ &  
	  0 & $[e^{2}]$& 0 & $[e^{3}]$&$0$ &$\cdots$ &0&$[e^d \wedge dr]$
	  \\[,1cm] \hline       &&&&&&&&&&&\\
	  $\per{p} \geq\per{k \! - \! 2} $ &  $1$ & 0 &  $[e]$ &  
	  0 & $[e^2]$ &  0 & $[e^{3}]$&0& $\cdots$&$[e^{d}]$& 0\\[,1cm]  \hline
	  \end{tabular}
	  \end{center}
      These calculations come directly from the equalities \refp{calculos}.

      We notice that the top dimensional basic cohomology group is 
	  isomorphic either to 0 or $\R$. These cohomology groups are finite 
	  dimensional. We recover the PoincarŽ Duality in the perverse sense:
      $$
      \lau{\IH}{*}{\per{p}}{\S^{k}/\mathcal{F}} \cong
      \lau{\IH}{n-*}{\per{q}}{\S^{k}/\mathcal{F}}
      $$
      for two complementary perversities: $\per{p} + \per{q} = 
      \per{t} =\per{k-3}$.
       
      This is more general. We have proved that the basic intersection cohomology is finite 
	  dimensional and verifies this perverse PoincarŽ Duality for the linear foliations 
	  (cf. \cite{SW2,SW3}) and for the riemannian flows (cf. \cite{RP}).

   \section{Top class}
      The top group $\lau{H}{n}{}{\mf}$ 
      of the basic cohomology of a  regular riemannian foliation $\F$ defined 
      on a connected compact 
  manifold $M$ is $\R$ or $0$. Here $n = 
      \codim_{M}\F$. We prove the same result for 
      the  top group $\lau{H}{n}{\per{p}}{\mf}$ when $\F$ is a SRF 
      defined on a connected compact manifold $M$.  To do so, we need three Lemmas.
      
      \bl
       \label{T1}
       Consider $(T_{S},\tau_{S},S)$ a foliated tubular neighborhood of a 
       singular
       stratum $S \in \SF$.  Fix $f \colon ]0,1[ \to [0,1]$ a 
       smooth function with with $f \equiv 0$ on 
 $]0,1/4]$ and $f\equiv 0$ on $[3/4,1[$. The map 
       $[\omega] \mapsto [df\wedge \omega]$ defines the isomorphism
       $$
\lau{\IH}{*-1}{\per{p},c}{D_{S}/\mathcal{F}} \cong
       \lau{\IH}{*}{\per{p},c}{D_{S} \times ]0,1[/\mathcal{F} \times 
       \mathcal{I}}.
       $$  
      \el
       \pro
       Consider the following differential complexes:
       \begin{itemize}
	\item[-] $\dos{\A}{*}  = 
	    \left\{ \om \in \lau{\Om}{*}{\per{p}}{D_S \times  ]0,3/4[/ \F\times 
	\mathcal{I}} \!  \Big/ \left[
	\begin{array}{l}
	\supp \om \subset K \times [c,3/4[\\
	\hbox{for a compact } K \subset D_S 
	\hbox{ and } 0 <c < 3/4 \! \! 
	\end{array} 
	\right.
	\right\}.
	$
	\item[-] $\dos{\B}{*} = 
	\left\{ \om \in \lau{\Om}{*}{\per{p}}{D_S \times  ]1/4,1[/ \F\times 
	\mathcal{I}}  \!  \Big/ \left[
	\begin{array}{l}
	\supp \om \subset K \times ]1/4,c]\\
	\hbox{for a compact } K \subset D_S 
	\hbox{ and } 1/4 < c <1 \! \! \! 
	\end{array} 
	\right.
	\right\}.
	$
	\item[-] $\dos{\C}{*} = 
	\left\{ \om \in \lau{\Om}{*}{\per{p}}{D_S \times ]1/4,3/4[/ \F\times 
	\mathcal{I}} \! \Big/  
	\left[\begin{array}{l}
	\supp \om \subset K \times ]1/4,3/4[\\
	\hbox{for a compact } K \subset D_S \! \! 
	\end{array}
	\right.
	\right\}.
	$
	\end{itemize}
       Proceeding as in \refp{MV} we get the short exact sequence
	$$
	0 \TO \lau{\Om}{*}{\per{p},c}{D_S \times ]0,1[/\F \times \mathcal{I}}
	\TO
	\dos{\A}{*}  \oplus \dos{\B}{*}
	\TO
	\dos{\C}{*}
	\TO 0.
	$$
	The associated long exact sequence is
       $$
       \cdots \to 
       \hiru{H}{i-1}{\dos{\C}{*}} 
       \stackrel{\delta}{\to}
       \lau{\IH}{i}{\per{p},c}{D_S \times ]0,1[/\F \times \mathcal{I}}
	\to 
	\hiru{H}{i}{\dos{\A}{*} } \oplus 
	 \hiru{H}{i}{\dos{\B}{*}}
       \to
	\hiru{H}{i}{\dos{\C}{*}}
       \to \cdots,
       $$
	where the connecting morphism is
	$
	\delta([\om]) = [df \wedge \om].
	$
       \medskip

       \smallskip
       
	   Before executing the calculation let us introduce some notation. 
       Let $\beta$ be a differential form on $\hiru{\Om}{i}{D_S \times 
       ]a,b[}$
       which does not 
       include the $dt$ factor. 
By $\Int{-}{c}\beta (s)\wedge ds$ and
       $\Int{c}{-}\beta  (s)\wedge ds$ we denote the forms on $\hiru{\Om}{i}{D_S \times 
       ]a,b[}$ 
       obtained from $\beta$ 
       by integration with respect to $s$, that is,
       $\left(\Int{-}{c}\beta (s)\wedge ds \right) 
       (x,t)(\vec{v}_{1}, 
       \ldots ,\vec{v}_{i}) = \Int{t}{c}(\beta (x,s)(\vec{v}_{1}, 
       \ldots ,\vec{v}_{i}) ) \ ds $
      and  on the other hand $\left( \Int{c}{-}\beta (s)\wedge ds 
       \right)(x,t)(\vec{v}_{1}, 
       \ldots ,\vec{v}_{i}) )= \Int{c}{t}(\beta 
       (x,s)(\vec{v}_{1}, 
       \ldots ,\vec{v}_{i}) ) \,  ds$
       where  $c \in ]a,b[$, $(x,t) \in D_S \times ]a,b[$ and 
       $(\vec{v}_{1}, 
       \ldots ,\vec{v}_{i}) \in T_{(x,t)}( D_{S} \times ]a,b[) $ .

       \smallskip

	Each differential form $\om \in \dos{\A}{*}$, $\dos{\B}{*}$ 
	and $\dos{\C}{*}$ 
	can be written $\om  = \alpha + dt \wedge \beta $ where $\alpha$ and 
	$\beta$ do not 
	contain $dt$. 
	Consider  a cycle $\om  = \alpha + dt \wedge \beta  \in \dos{\A}{i}$ with $\supp \om \subset K \times [c,3/4[$ for a 
	compact 
	$K \subset D_{S}$ and $0 < c < 3/4$.
	We have 
	$
	\om = 
	-d \left( \Int{-}{c/2} \beta (s) \wedge  ds\right).
	$
	Since $ \supp \Int{-}{c/2} \beta (s) 
	\wedge ds \subset K \times [c,3/4[$ then we get
	$
	\hiru{H}{*}{\dos{\A}{\cdot} } =0$. In the same way, we 
	get  $
	\hiru{H}{*}{\dos{\B}{\cdot} } =0$. 
We conclude that $$\delta  \colon \hiru{H}{*-1}{\dos{\C}{\cdot} } 
\TO
       \lau{\IH}{*}{\per{p},c}{D_S \times ]0,1[/\F \times \mathcal{I}}
$$ is an isomorphism.

\smallskip 
The foliated homotopy $L \colon D_{S} \times ]1/4,3/4[ \times [0,1] \to D_{S} \times 
]1/4,3/4[ $,
defined by $L((u,t),s) = (u,(1-s)t+s/2)$,  verifies $L(K \times 
]1/4,3/4[ \times [0,1]) \subset K \times 
]1/4,3/4[ $ for each compact $K \subset D_{S}$.
So, we get that the map $[\alpha] \mapsto [\alpha]$ induces the 
isomorphism

\begin{equation}
    \label{I}
I \colon \lau{\IH}{*}{\per{p},c}{D_{S}/\mathcal{F}} \TO \hiru{H}{*}{\dos{\C}{\cdot} } .
\end{equation}

\nt The isomorphism $\delta \rondp I$ gives the result.
	\qed
 \bl
 \label{T2}
 Consider $T_S$ a foliated tubular neighborhood of a 
 singular stratum $S \in \SF$.  The inclusion $T_S\backslash S \hookrightarrow  T_S$ induces the onto map
 $$
 \iota \colon\lau{\IH}{n}{\per{p},c}{(T_S\backslash  S)/\mathcal{F}} \to
 \lau{\IH}{n}{\per{p},c}{T_S/\mathcal{F}}.
 $$
 Moreover, if $\per{p}(S) \leq \per{t}(S)$ then $\iota$ is an isomorphism.
\el
 \pro
We proceed in several  steps.

 \Zati {\em Rewriting $\lau{\IH}{*}{\per{p},c}{T_S/\mathcal{F}}$}. 
 
 \smallskip
 
 \nt We have seen in 
 \cite{SW2} that we can identify $\lau{\Om}{*}{\per{p},c}{T_S/\mathcal{F}}$ with the 
 complex
 $$
 \left\{ \om \in \lau{\Om}{*}{\per{p},c}{D_S \times [0,1[ / \mathcal{F} \times 
 \mathcal{I}} \ \big|  \ \|\om |_{D_S \times \{0 \}} \|_{\tau_S}\leq \per{p}(S) \hbox{ and } 
 \|d\om |_{D_S \times \{0 \}}\|_{\tau_S} \leq \per{p}(S)\right\}.
 $$
  Here $\| - \|_{\tau_S}$ denotes the 
 {\em vertical degree}  relatively to the fibration $\tau_S \colon D_S \equiv D_S 
 \times \{ 0 \} \to S$. Recall that $\| 0\|_{\tau_S} = -\infty$.
Under this identification, the complex 
$\lau{\Om}{*}{\per{p},c}{(T_S\backslash  S)/\mathcal{F}}$ 
becomes
$
\lau{\Om}{*}{\per{p},c}{D_{S} \times ]0,1[, \mathcal{F} \times 
\mathcal{I}}.
$

 \zati {\em Chasing $\lau{\IH}{n}{\per{p},c}{T_S/\mathcal{F}}$}. 
 
 \smallskip
 
 \nt Consider the complex
  \begin{itemize}
     \item[-] 
     $
     \dos{\D}{*}  = 
     \left\{ \om \in \lau{\Om}{*}{\per{p}}{D_S \times [0,3/4[ / \mathcal{F} \times 
 \mathcal{I}} \ \Big/ \ \left[
 \begin{array}{l}
 \|\om |_{D_S \times \{0 \}} \|_{\tau_S}\leq \per{p}(S), \\[,1cm]
 \|d\om |_{D_S \times \{0 \}}\|_{\tau_S} \leq \per{p}(S) \hbox{ and } \\[,1cm]
 \supp \om \subset K \times [0,3/4[ \\
 \hbox{for a compact } K \subset D_S
 \end{array} 
 \right.
 \right\}
 $
  \end{itemize}
 and the complexes
$\dos{\B}{*}$, $\dos{\C}{*}$ as in the proof 
of Lemma \ref{T1}.
 The short exact sequence
 $$
 0 \TO \lau{\Om}{*}{\per{p},c}{T_S/\F}  
 \TO
 \dos{\D}{*}  \oplus \dos{\B}{*}
 \TO
 \dos{\C}{*}
 \TO 0
 $$
 produces the long exact sequence
 $$
 \to
 \hiru{H}{n-1}{\dos{\D}{*} } \oplus \hiru{H}{n-1}{\dos{\B}{*}}
\to
 \hiru{H}{n-1}{\dos{\C}{*}}
 \stackrel{\delta'}{\TO}
 \lau{\IH}{n}{\per{p},c}{T_S/\mathcal{F}}
 \to
 \hiru{H}{n}{\dos{\D}{*} } \oplus \hiru{H}{n}{\dos{\B}{*}}
\to
 \hiru{H}{n}{\dos{\C}{*}}
 \to
 0
 $$
 where the connecting morphism is
 $
 \delta'([\om ]Ê) = [df\wedge \om]
 $
 for a smooth map $f\colon [0,1[ \to [0,1]$ with $f \equiv 0$ on 
 $[0,1/4]$ and $f\equiv 0$ on $[3/4,1[$.

 \zati {\em Relating
 $\hiru{H}{*}{\dos{\C}{\cdot}}$ and 
 $\hiru{H}{*}{\dos{\D}{\cdot}}$}. 
 
 \smallskip
 
\nt We have already seen that  $
\hiru{H}{*}{\dos{\B}{\cdot}}
= 0$.  Consider  a cycle $\om  = \alpha + dt \wedge \beta \in \dos{\D}{n}$. For degree reasons, we have $\alpha(0) = 0$ and then
 $
 \om = 
 d \left( \Int{0}{-} \beta (s) \wedge  ds\right).
 $
 Since $ \supp \Int{0}{-}  \beta (s) 
 \wedge ds \subset K \times [0,3/4[$ we get 
 $\hiru{H}{n}{\hiru{D}{\cdot}{ [0,3/4[ }} =0$.
 From the above long exact sequence, we obtain the exact sequence
 $$
  \hiru{H}{n-1}{\dos{\D}{*} } 
 \TO
  \hiru{H}{n-1}{\dos{\C}{*}}
  \stackrel{\delta'}{\TO}
  \lau{\IH}{n}{\per{p},c}{T_S/\mathcal{F}}
  \TO
  0.
  $$

 \zati {\em Conclusion}.
 
 \smallskip
 
 \nt Since the map $I$ is an isomorphism (cf. \refp{I}) the 
 composition
 $
\delta'  \rondp I\colon 
 \lau{\IH}{*-1}{\per{p},c}{D_{S}/\F} \to 
 \lau{\IH}{n}{\per{p},c}{T_S/\mathcal{F}}
 $
 is an onto map. The Lemma \ref{T1} gives that $\iota  $ is an onto map. 
It remains to prove that the map $\iota$ is an isomorphism when $\per{p}(S) \leq 
 \per{t}(S)$. We prove   $
  \hiru{H}{n-1}{\dos{\D}{\cdot}} =0$.
  Proceeding as in $(c)$ it suffices to consider a cycle
$\om  = \alpha + dt \wedge \beta \in \hiru{D}{n-1}{ 
  [0,3/4[}$ and show that $\alpha(0) = 0$.
  
  Let us suppose that $\alpha (0) \not= 0$. 
  The  foliation $\F$ induces a foliation  $\F_{\tau_{S}}$ tangent to the 
  fibers of $\tau_{S} \colon D_{S }\to S$  such that $\dim \F = \dim 
  \F_{\tau_{S}} + \dim \F_{S}$
  (cf. \cite{Mo,SW1}). By degree reasons, 
  since $\alpha (0) \in \lau{\Om}{n-1}{\per{p},c}{D_{S}/\F}$, we can 
  write
  $$
  \per{t}(S) \geq \per{p}(S) \geq ||\alpha (0)||_{\tau_{S}} =  \underbrace{(\dim M - \dim 
  S) -1}_{dimension \ of  \ the  \ fiber \  of  \ \tau_{S}}  - 
  \underbrace{(\dim \F - \dim \F_{S})}_{dim \ \F_{\tau_{S}}}= 
  \per{t}(S) +1 .
  $$
  This contradiction gives $\alpha (0) =0$.
 \qed
 
 \bl
 \label{T3}
Let $N$ be a compact manifold endowed with a RF $\mathcal{N}$.
If $O$ is a connected saturated open subset of $N$, then 
 $
 \lau{H}{n}{c}{O/\mathcal{N}} = 0 \hbox{ or }\R .
 $
 \el
 \pro
We proceed in two steps.
 
 \medskip
 
 \underline{\em The foliation $\mathcal{N}$ is transversally orientable}. 
 The result comes essentially from the basic Poin\-ca\-rŽ duality theorem for 
 non-compact manifolds (cf.\cite{Se}). 
 The foliation $\mathcal{N}_{O}$ is a transversally 
 orientable complete foliation since $O$ is a saturated open subset of 
 the compact manifold $N$. Then we have the isomorphism 
 $\lau{H}{n}{c}{O/\mathcal{N}}  \cong 
 \hiru{H}{0}{O/\mathcal{N},\mathcal{P}}$ where $\mathcal{P}$ is the homological 
 orientation sheaf of $\mathcal{N}_{O}$.
 Since the manifold $O$ is connected, the  sheaf $\mathcal{P}$ is 
 locally trivial and the stalk is $\R$ then $\lau{H}{n}{c}{O/\mathcal{N}} = 0 \hbox{ or }\R .$

 \medskip
 
 \underline{\em General case}. 
 Consider the 
  transverse orientation covering $\pi \colon (N^*,\mathcal{N}^*) \to 
  (N,\mathcal{N})$ (see \cite{HH}). The covering is given by a foliated 
  action of $\Z_{2}$. The foliation $\mathcal{N}^*$ is 
  a transversally orientable RF. We have the equality 
  $
  \lau{H}{n}{c}{O/\mathcal{F}} = \left( \lau{H}{n}{c}{\pi^{-1}(O)/\mathcal{F}^*}\right)^{\Z_{2}}.
  $
  The subset $\pi^{-1}(O)$ 
  is a saturated open subset of $N^{*}$. 
  If $\pi^{-1}(O)$ is connected then the result comes from the 
  previous case.
  If $\pi^{-1}(O)$ is not connected then $\pi^{-1}(O)$ has two 
  connected components foliated diffeomorphics to $O$ and  the 
  $\Z_{2}$-action interchanges them. The result comes now from the 
  previous case.
  
 \qed
 
 The first  result of this Note is the following.

 \bt
 \label{tteeoo}
 Let $M$ be a connected compact manifold endowed with a SRF $\F$. 
 If $n = \codim_{M} \F$ and $\per{p}$ a perversity on $M$,
 then
 $$
 \lau{\IH}{n}{\per{p}}{\mf} = 0 \hbox{ or } \R,
 $$
 \et
 \pro
For each $i\in \Z$ we 
 write: 
 \begin{itemize}
     \item[-] $\Sigma_{i} \subset M$ the union of strata whose dimension is 
 less or equal than $i$,
 \item[-] $T_{i}$ the tubular neighborhood of $\Sigma_{i}$ in $M 
 \backslash  \Sigma_{i-1}$. 
 
 \end{itemize}
 We have $M \backslash   \Sigma_{-1} = M$ and $M \backslash   \Sigma_{m-1} = R$, 
 where $m = 
 \dim M$. The Lemma \ref{T3} gives $\lau{H}{n}{c}{R/\F} =0 \hbox{ or } \R$.
 We get the result if we prove that the inclusion 
 $
 M\backslash   \Sigma_{i}
 \hookrightarrow
 M\backslash   \Sigma_{i-1}
 $
 induces an onto map
 $\lau{H}{n}{c}{(M\backslash   \Sigma_{i})/\F} 
 \to  \lau{H}{n}{c}{(M\backslash   
 \Sigma_{i-1})/\F}$, for $i \in \{ 0, \ldots , m-1\}$.

 From the open covering
 $
 {\displaystyle
 \left\{ M\backslash  \Sigma_i, \bigcup_{\dim S =i}T_{S}\right\}
 }
 $
 of $M\backslash  \Sigma_{i-1}$, 
 we obtain the Mayer-Vietoris sequence
 $$
 \bigoplus_{\dim S =i} \lau{H}{n}{c}{(T_{S}\backslash  S)/\F} \TO
 \lau{H}{n}{c}{(M\backslash  \Sigma_i)/\F}
 \oplus
 \bigoplus_{\dim S =i} \lau{H}{n}{c}{T_{S}/\F} 
 \TO
 \lau{H}{n}{c}{(M\backslash  \Sigma_{i-1})/\F} \to 0.
 $$
 Now, the Lemma \ref{T2} gives the result.
 \qed

Combining Theorem \ref{tteeoo} with the tautness characterization 
of \cite{RSW1}, we get the following Corollary.

 \bc
 \label{T}
 Let $M$ be a connected compact manifold endowed with a SRF $\F$. Let 
 us suppose that $\F$ is transversally orientable. 
Consider  $\per{p}$ a perversity on $M$ with 
  $\per{p}\leq \per{t}$.
If $n = \codim_{M} \F$, then the two following statements are equivalent:

 \Zati The foliation $\F_{R}$ is taut, where $R$ is the regular 
 stratum of $(M,\F)$.

 \zati The cohomology group $\lau{\IH}{n}{\per{p}}{M/\F}$ is $\R $.

 \bigskip

 \ec
    \pro
   We know from \cite{RSW1}  that the condition $(a)$ is equivalent to 
   $\lau{H}{n}{c}{R/\F} = \R $. So, it suffices to prove 
   that $\lau{H}{n}{c}{R/\F} \cong \lau{\IH}{n}{\per{p}}{M/\F}$ (cf. 
   Lemma \ref{T3}).  We proceed as in the proof of the previous Theorem 
   changing ``onto map" by ``isomorphism".
   \qed

   \prg {\bf Remarks.}

   \Zati The perversity  $\per{p} = \per{-\infty}$ verifies $\per{p}\leq 
   \per{t}$. In this case the group
   $\lau{\IH}{n}{\per{p}}{\mf}$ becomes $\lau{H}{n}{}{\mf, 
   \Sigma/\F}$. Here, the  relative basic cohomology  
   $\lau{\IH}{*}{\per{p}}{\mf, \Sigma/\F}$
   is computed from the relative basic 
   complex 
   $
   \lau{\Om}{*}{}{\mf, \Sigma/\F} = \{ \om \in \lau{\Om}{*}{}{\mf} 
   \  | \ \om \equiv 0 \hbox{ on } \Sigma\}.
   $

   \zati The boundary perversity  $\per{\partial}$ verifies $\per{\partial}\leq 
   \per{t}$. In this case the group
   $\lau{\IH}{n}{\per{\partial}}{\mf}$ becomes $\lau{H}{n}{}{\mf, 
    \partial(M/\F)}$. Here, the  relative basic cohomology  
    $\lau{H}{*}{}{\mf,\partial(M/\F)}$
   is computed from the relative basic 
   complex 
   $
   \lau{\Om}{*}{}{\mf,\partial(M/\F)} = \{ \om \in \lau{\Om}{*}{}{\mf} 
   \  | \ \om \equiv 0 \hbox{ on the boundary strata}\}.
   $
   In particular, when the boundary strata do not appear then
   $\lau{\IH}{n}{\per{\partial}}{\mf} = \lau{H}{n}{}{\mf}$.

   \zati When $\per{p} \not\leq \per{t}$ then the group 
   $\lau{\IH}{n}{\per{p}}{\mf} $ does not establish the tautness 
   of $(R,\F)$. For example, we always have 
  $\lau{\IH}{n}{\per{t} + \per{1}}{\mf}  = \hiru{H}{n}{(M\backslash  
  \Sigma)/\F} =0$ if $\Sigma \neq\emptyset$.
    
   \zati In the Theorem \ref{tteeoo} and the Corollary \ref{T} we can suppose that the manifold $M$ 
   is not compact and replace $\lau{\IH}{n}{\per{p}}{\mf}$ by 
   $\lau{\IH}{n}{\per{p},c}{\mf}$.

\end{document}